\documentclass[reqno]{amsart}
\usepackage{amsthm,wrapfig,amsmath,amssymb,amsfonts,setspace,verbatim,graphicx,mathtools,mathrsfs,commath,float,mdframed,frame,xcolor,qtree,enumitem, ulem,afterpage}

\usepackage[hidelinks]{hyperref}
\usepackage[T1]{fontenc}

\DeclareMathOperator{\diam}{diam}

\usepackage{newunicodechar}
\usepackage{lmodern}

\usepackage[margin=1.4in]{geometry}

\newtheorem{ut}{Theorem}

\newtheorem{uc}[ut]{Corollary}

\begin{document}

\title[Totally disconnected subsets of  chainable continua]{Totally disconnected subsets of  chainable continua}

\address{Department of Mathematics, Auburn University at Montgomery}
\email{dsl0003@auburn.edu}
\author[D.S. Lipham]{David S. Lipham}

\begin{abstract}We show that the endpoint set of a Suslinian chainable continuum must be zero-dimensional at some point. In particular, it cannot be homeomorphic to complete Erd\H{o}s space. This answers a question of Jerzy  Krzempek.\end{abstract}

\maketitle


A \textbf{continuum} is a compact connected metric space.  A continuum $X$ is  \textbf{Suslinian} if there is no uncountable collection of pairwise disjoint non-degenerate subcontinua of $X$. A continuum $X$ is \textbf{chainable} if for every $\varepsilon>0$ there are finitely many open sets $U_1,\ldots ,U_n$ covering $X$ such that $\diam (U_i)<\varepsilon$  and $U_i\cap U_j\neq \varnothing$$ \iff$$|i-j| \leq 1$ for all $i,j\leq n$.  A point $x$ in a chainable continuum $X$ is an \textbf{endpoint} of $X$ if the $\varepsilon$-chains $ U_1,\ldots ,U_n$ can always be chosen so that $x\in U_1$. For equivalent definitions of endpoints see \cite[Theorem 13]{bing} or \cite[p.609]{aki}. The set of all endpoints of $X$ is denoted $E(X)$.

A topological space $X$ is:
\begin{itemize}\renewcommand{\labelitemi}{\scalebox{.5}{{$\blacksquare$}}}
\item  \textbf{totally disconnected} if $X$ does not have any  connected subset with more than one point; 
\item  \textbf{zero-dimensional at} $x\in X$ if  every neighborhood of $x$ contains a clopen neighborhood of $x$; 

\item \textbf{zero-dimensional} if $X$ is zero-dimensional at each of its points (i.e.\ $X$ has a basis of clopen sets). \end{itemize}

\renewcommand*{\thefootnote}{$\dagger$}

Jerzy Krzempek has shown  that for every  zero-dimensional Polish space $E$   there is a  Suslinian chainable continuum $X$ such that $E(X)$ is homeomorphic to $E$ \cite[Theorem 1]{krz}.  So interestingly enough there is a Suslinian chainable continuum whose endpoint set is homeomorphic to the space of irrational numbers. Conversely,  if $X$ is Suslinian and chainable then $E(X)$ is totally disconnected \cite[Theorem 11]{aki}\footnote{The theorem referenced here states only that $E(X)$ does not contain any non-degenerate continuum, but its proof can be easily modified to show that $E(X)$ does not contain any non-degenerate connected set. In fact, supposing that $C\subset E(X)$ is connected and non-degenerate, a contradiction can be reached using the fact that the continuum $\overline C$ is decomposable.} and Polish \cite[Theorem 4.2]{do}.   It is unknown whether $E(X)$  must be zero-dimensional  \cite[Problem 1]{krz}. Krzempek in particular asked whether $E(X)$ could be homeomorphic to complete Erd\H{o}s space $\mathfrak E_{\mathrm{c}}$, a famous example of a totally disconnected Polish space which is not zero-dimensional at any point \cite{erd,dij}. In this paper we provide a negative answer with the following.

\begin{ut}If $X$ is a Suslinian chainable continuum, then every totally disconnected $G_{\delta}$-subset of $X$ is zero-dimensional at some point.\end{ut}
 
\begin{uc} If $X$ is a Suslinian chainable continuum, then $E(X)$ is zero-dimensional at some point. \end{uc}
  
\begin{uc}Every  chainable continuum that  contains  $\mathfrak E_{\mathrm{c}}$ is non-Suslinian.  In particular, $\mathfrak E_{\mathrm{c}}$ is not homeomorphic to the endpoint set of a Suslinian chainable continuum. \end{uc}

\renewcommand*{\thefootnote}{$\ddagger$}

\subsection*{Examples}The following examples show that Theorem 1 and Corollary 3 essentially cannot be improved upon. 
\begin{itemize}\renewcommand{\labelitemi}{\scalebox{.5}{{$\blacksquare$}}}
\item[(A)] There is a Suslinian chainable continuum that contains a totally disconnected  Polish space of positive dimension. This example was constructed by Howard Cook and Andrew Lelek  \cite[Example 4.2] {iv}.  By \cite[Theorem 3.2] {iv} it contains a totally disconnected set $P$ with countable complement, and by the statement of \cite[Example 4.2] {iv} the set $P$ is not zero-dimensional.\footnote{Note that in \cite{iv}, the term \textit{hereditarily disconnected} is used instead of totally disconnected, and \textit{totally disconnected} is defined to be a stronger condition which is still satisfied by  every zero-dimensional space.} 
\item[(B)] There is a hereditarily decomposable chainable continuum which homeomorphically contains  $\mathfrak E_{\mathrm{c}}$  (e.g.\  the Cantor organ \cite[p.\;191]{kur}).
\item[(C)] There is  a Suslinian dendroid whose endpoint set is homeomorphic to  $\mathfrak E_{\mathrm{c}}$. This example can be obtained as a quotient of the Lelek fan and is due to  Piotr Minc and Ed Tymchatyn (personal communication). 
\end{itemize}
 \begin{figure}[h]\includegraphics[scale=0.042]{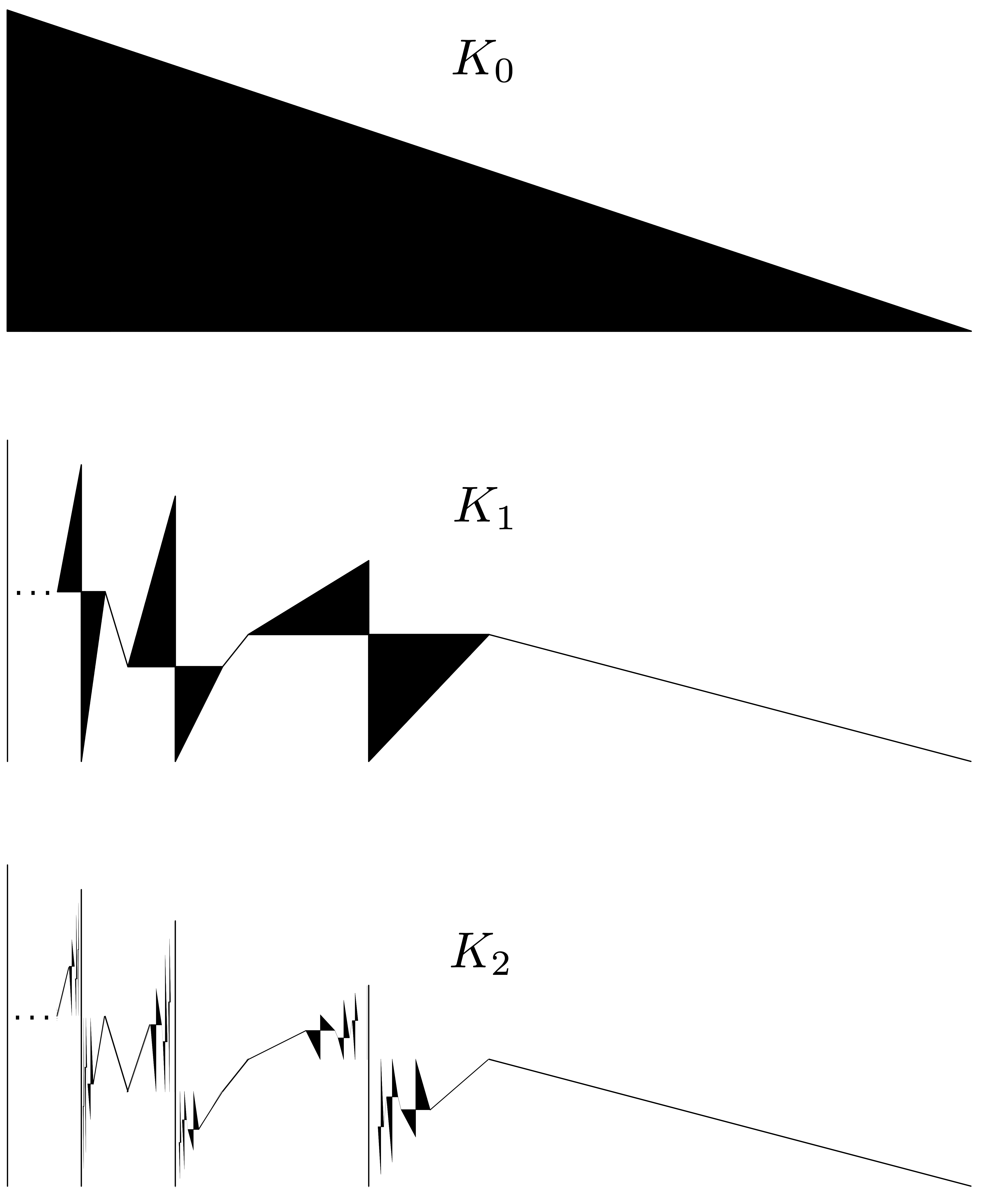}
 \caption{Illustration of Example A. The set $K_1$ is constructed so that the points where two triangles intersect accumulate onto the middle third of the left-most segment $S$. This principle is repeated in the smaller triangles to form $K_2, K_3$, etc.  The set $X=\bigcap _{n=0}^\infty K_n$ is a Suslinian chainable continuum.  If  $Q$ is any countable set which intersects each arc of $X$, then $P=X\setminus Q$ is totally disconnected. But every clopen subset of $P$ meeting $S$  must contain all of $P\cap S$ (proved by Cook and Lelek). Hence $P$ is not zero-dimensional.}
\end{figure}

\subsection*{Proof of Theorem 1} Let  $X$ be a Suslinian chainable continuum.  Let $E$  be a totally disconnected $G_{\delta}$-subset of $X$. We will produce a point $e$ at which $E$ is zero-dimensional. To that end, let $\mathcal K$ be the set of all non-degenerate connected components of $\overline E$. Note that $\mathcal K$ is countable by  the Suslinian property. So by Baire's theorem either $\bigcup \mathcal K$ has empty interior in $\overline E$, or there exists $K\in\bigcup \mathcal K$ which contains a non-empty open subset of $\overline E$.
 \medskip
 
 \textbf{Case 1:} $\bigcup \mathcal K$ has empty interior in $\overline E$. 
 \medskip
 
 Then $\overline E\setminus \bigcup \mathcal K$ is a dense $G_{\delta}$-set in $\overline E$. Since $E$ is also dense  and $G_{\delta}$ in $\overline E$,  there exists $e\in E\cap (\overline E\setminus \bigcup \mathcal K)$.  Then  $\{e\}$ is a connected component of the compactum $\overline E$. So $\overline E$ is zero-dimensional at $e$ (cf.\ \cite[Section 1.4]{eng}). Therefore $E$ is zero-dimensional at $e$.
  \medskip
  
 \textbf{Case 2:} There exists $K\in\bigcup \mathcal K$ which contains a non-empty open subset of $\overline E$.
 \medskip
  
Let $U$ be a non-empty relatively open subset of $\overline E$ that is contained in $K$. 

The continuum $K$ is Suslinian and chainable, as these properties are inherited from $X$. Therefore $K$ is hereditarily decomposable \cite[Theorem 1.1]{lel} and irreducible \cite[Theorem 12.5]{nad}. By Kuratowski's  theory of tranches (see  \cite[\S 48]{kur} or \cite[p.15]{es}), there is a mapping $\varphi:K\to [0,1]$ such that $\varphi^{-1}\{t\}$ is a nowhere dense subcontinuum of $K$ for every $t\in [0,1]$. By the Suslinian property, the set $$Q=\big\{t\in [0,1]:|\varphi^{-1}\{t\}|>1\big\}$$  is countable. Therefore   $K\setminus \varphi^{-1}(Q)$ is dense $G_{\delta}$ in $K$.  So there exists $$e\in  E\cap U\setminus \varphi^{-1}(Q).$$ 
By  compactness of $K$  there exist $a,b\in [0,1]\setminus Q$ such that  $$a<\varphi(e)<b$$ and $\varphi^{-1}[a,b]\subset U$.   

We claim that there is a relatively  clopen  $C\subset E\cap \varphi^{-1}[a,\varphi(e)]$ which contains $e$ and misses the point $\varphi^{-1}(a)$. To see this, note that $E\cap \varphi^{-1}[a,\varphi(e)]$ is non-degenerate because $E$ is dense in $\varphi^{-1}(a,\varphi(e))$. Since $E$ is totally disconnected,  $E\cap \varphi^{-1}[a,\varphi(e)]$ is not connected. Let $W$ be a proper clopen subset of the space $E\cap \varphi^{-1}[a,\varphi(e)]$ such that $e\in W$. Since $K\setminus \varphi^{-1}(Q)$   is dense in $\varphi^{-1} (a,\varphi(e))$, and $W$ is a proper closed subset of $E\cap \varphi^{-1}[a,\varphi(e)]$,   there exists $c\in  [a,\varphi(e))\setminus Q$ such that $\varphi^{-1}(c)\notin W$. Then $C=W\cap \varphi^{-1}[c,\varphi(e)]$ is as desired.

Likewise there is a relatively clopen $D\subset E\cap \varphi^{-1}[\varphi(e),b]$ which contains $e$ and misses the point $\varphi^{-1}(b)$. Then $C\cup D$ is a clopen subset of $E$ that contains $e$ and lies inside of $U$. If $U'$ is any smaller open subset of $\overline E$ containing $e$, then the same argument will produce a clopen subset of $E$ which contains $e$ and lies inside of $U'$. Therefore $E$ is zero-dimensional at $e$. 

This concludes the proof of Theorem 1. \hfill $\blacksquare$

\subsection*{Remarks} Corollaries 2 and 3 follow from Theorem 1 because $E(X)$ and $\mathfrak E_{\mathrm{c}}$ are totally disconnected Polish spaces, and $\mathfrak E_{\mathrm{c}}$ is nowhere zero-dimensional. 
 
By applying Theorem 1 locally, it can be seen that if $E$ is any totally disconnected $G_{\delta}$-subset of $X$ (a Suslinian chainable continuum), then the set of points at which $E$ is zero-dimensional is dense in $E$. The set of points at which a separable metrizable space is zero-dimensional is always $G_{\delta}$, so we get the following strengthening of Corollary 2: \textit{The set of points at which $E(X)$ is zero-dimensional is dense $G_{\delta}$ in $E$.} Thus $E(X)$ is zero-dimensional ``almost everywhere''.






\end{document}